\documentclass{amsart}
\usepackage{amssymb,amsmath}

\newcommand{\RR}{{\mathbb R}}
\newcommand{\NN}{{\mathbb N}}
\newcommand{\CC}{{\mathbb C}}
\newcommand{\C}{{\mathbb S}}
\newcommand{\DD}{{\mathbb D}}

\newcommand{\e}{\varepsilon}

\newcommand{\del}{\partial}
\newcommand{\om}{\omega}
\newcommand{\Om}{\Omega}

\newcommand{\si}{\sigma}

\newcommand{{\loc}}{{\ell\mathrm oc}}

\newtheorem{theorem}{Theorem}
\newtheorem{lemma}{Lemma}
\newtheorem{proposition}{Proposition}
\newtheorem{corollary}{Corollary}

\newtheorem{definition}{Definition}

\newtheorem*{main}{Main Theorem}

\begin{document}
\title[Stability for Dini Continuous Conductivities]
{Stability of the Inverse Problem for Dini Continuous Conductivities in the Plane}

\author{Robert McOwen }

\author{Bindu K Veetel}

\date{March 15, 2020}

\maketitle

\begin{abstract}  
We show that the inverse problem of Calderon for conductivities in a two-dimensional Lipschitz domain is stable in a class of conductivities that are Dini continuous. This extends previous stability results when the conductivities are known to be H\"older continuous.
\end{abstract}

\addtocounter{section}{0}
\addtocounter{lemma}{0}

\section{Introduction}

\noindent
We  study the inverse problem of Calderon in two-dimensions, namely the determination of the conductivity function $\gamma$ in a bounded domain $U\subset\RR^2$ from the Dirichlet to Neumann (or ``voltage to current'') map 
\begin{equation}\label{def:DtNmap}
\Lambda_\gamma: H^{1/2}(\partial U)\to H^{-1/2}(\partial U),
\end{equation}
which associates to the Dirichlet data $f\in H^{1/2}(\partial U)$ the function 
$\Lambda_\gamma(f)=\gamma\frac{\partial u}{\partial\nu}$
where $\nu$ is the exterior unit normal on $\partial U$ and $u$ is the unique solution of the Dirichlet problem
\begin{equation}\label{eq:conductivity}
\begin{aligned}
\nabla \cdot \gamma \nabla u &= 0 \quad\hbox{in}\ U, \\
u&=f \quad \hbox{on}\ \partial U.
\end{aligned}
\end{equation}
The uniqueness of the inverse of $\gamma\mapsto \Lambda_\gamma$ for conductivities $\gamma\in L^\infty(U)$ was proved in \cite{AP}: this means that $\Lambda_{\gamma_1}=\Lambda_{\gamma_2}$ implies $\gamma_1=\gamma_2$. Stability is the stronger statement that
\begin{equation}\label{def:stability}
\|\gamma_1-\gamma_2\|_{L^\infty}\leq V(\|\Lambda_{\gamma_1}-\Lambda_{\gamma_2}\|_{\partial U})
\end{equation}
where $\|\cdot\|_{\partial U}$ denotes the operator norm $H^{1/2}(\partial U)\to H^{-1/2}(\partial U)$ and $V(\rho) > 0$ is a non-decreasing ``stability function'' that satisfies $V(\rho)\to 0$ as $\rho\to 0$. It is known that stability (\ref{def:stability}) does not hold for all $\gamma\in L^\infty(U)$ (see the example of Alessandrini \cite{A}), so it becomes interesting to inquire what is the most general class of functions $\gamma$ for which stability holds?

In \cite{BFR} it was shown that if $\partial U$ is Lipschitz and the $\gamma_i$ are H\"older continuous functions on $\overline U$ satisfying $\e <\gamma_i <1/\e$, then (\ref{def:stability}) holds with $V(\rho)=C\,(\log\rho)^{-a}$ for some $a>0$. In this paper we want to prove stability for Dini continuous conductivities. Because we need fairly sharp estimates, we will assume that our conductivities have modulus of continuity $\varpi$  of the form
\begin{equation}\label{eq:bar-om}
\varpi(r)=|\log r|^{-\alpha} \quad\hbox{for}\ 0<r\leq 1/2.
\end{equation}
\noindent
Note that $\varpi$ satisfies the Dini condition $\int_0^{1/2} r^{-1}\varpi(r)\,dr<\infty$ if $\alpha>1$; 
however, for technical reasons, we shall require $\alpha>3/2$.

\begin{main}
Suppose $U$ is a bounded, Lipschitz domain and $\gamma_1,\gamma_2$ are functions on $\overline{U}$ which have modulus of continuity  $\varpi$ as in \eqref{eq:bar-om} with $\alpha>3/2$
and satisfy $\e <\gamma_i(z) <1/\e$ for some $\e>0$ and all $z\in U$. Then there exists a stability function $V(\rho)$ satisfying $V(\rho)\to 0$ as $\rho\to 0$ such that (\ref{def:stability}) holds.
\end{main}
\noindent
The proof of this theorem builds on the ideas of  \cite{AP} and  \cite{BFR}. However, to handle functions with modulus of continuity \eqref{eq:bar-om}, we require some modified Bessel potentials that we have not seen  in the literature. As usual, Bessel potentials for $p\not=2$ are somewhat tricky, so we here consider just $p=2$; this leads to the restriction $\alpha>3/2$.


\section{Some Function Spaces}

Recall that a {\it modulus of continuity} $\om(r)$ is a function $\om:[0,\infty)\to[0,\infty)$ which is strictly increasing for $r$ near $0$ and satisfies $\om(0)=0$. We shall also assume
\begin{equation}
\om(r)\geq c\,r^\e \ \hbox{ for some $\e\in (0,1)$ and all $0<r\leq 1/2$,}
\end{equation}
and, for convenience, we assume $\om$ is constant for $r\geq 1/2$.
For any bounded domain $U$ in $\RR^2$, we  define the following function spaces: 
\begin{definition} Let $C^\om(\overline U)$ denote the Banach space of functions $f\in C(\overline U)$ for which
$|f(x)-f(y)|\leq C\,\om(|x-y|)$ for all $x,y\in \overline U$ with the norm 
\begin{equation}\label{defnorm}
\|f\|_{C^\om(\overline U)}:=\sup_{x\in \overline U}|f(x)|+ \sup_{\substack{ x,y\in \overline U \\ x\not= y}}\frac{|f(x)-f(y)|}{\om(|x-y|)} .
\end{equation} 
 Let $C^{1,\om}(\overline U)$ denote the Banach space of functions $f\in C^1(\overline U)$ whose first order derivatives $\del f/\del z $ and $\del f/ \del \overline z $ are in $C^\om(\overline U)$ with the norm
\begin{equation}\label{defnorm1}
\|f\|_{C^{1,\om}(\overline U)}:=\|\del f\|_{C^\om(\overline U)} + \|\overline \del f\|_{C^\om(\overline U)} + \| f\|_{C^0(\overline U)}.
\end{equation}
\end{definition}
\begin{definition}
For any domain $\Om$, let  $C^\om(\Om)$ and $C^{1,\om}(\Om)$ denote respectively the union of all $C^\om(\overline U)$ and $C^{1,\om}(\overline U)$ where $U$ is compactly contained in $\Om$.
We also let $C_0^\om(\Om)$ and $C_0^{1,\om}(\Om)$ denote those functions with compact support in $\Omega$. Frequently we take $\Omega=\DD$, the unit disk, or $\DD_R:=\{x:|x|< R\}$.
\end{definition}
\noindent
When $\omega(r)=r^\gamma$ for $\gamma\in(0,1)$, then $C^\omega(\Omega)$ is traditionally written as $C^\gamma(\Omega)$, the functions which are {\it H\"older continuous of order $\gamma$}. For another example, we can extend the function $\varpi(r)$ as in \eqref{eq:bar-om} to be constant on $[1/2,\infty)$ and, since $\varpi(r)\to 0$ as $r\to 0$, we can define $\varpi(0)=0$ to make $\varpi(r)$ a modulus of continuity. Thus we may consider the function spaces $C^\varpi(\Omega)$, etc. Notice that $C^\varpi(\Omega)$ is larger than any H\"older space $C^\gamma(\Omega)$ for $\gamma\in (0,1)$.

\medskip
For $1\leq p<\infty$, we let $H^{1,p}(\Omega)$ denote the 1st-order $L^p$-Sobolev space for $\Omega$
and $H_\loc^{1,p}(\Omega)$  functions that are in $H^{1,p}({U})$ for any compact subset $\overline{U}\subset \Omega$. However, we will also be interested in less regular functions.
Suppose $\vartheta:[0, \infty) \to [0, \infty)$ is  increasing with $\vartheta(r)\to\infty$ as $r\to\infty$. 
Then we use the  Fourier transform $\widehat f(\xi)=\int_{\RR^2}e^{2\pi ix\cdot x}f(x)\,dx$ to define the following Banach (and Hilbert) space.
\begin{definition} 
For  $\vartheta:[0, \infty) \to [0, \infty)$ nondecreasing with $\vartheta(r)\to\infty$ as $r\to\infty$,
\begin{equation}\label{def:W^(Omega,2)}
W^{\vartheta,2}(\RR^2):=\{ f\in L^2(\RR^2): \|f\|_{W^{\vartheta,2}}^2:=\int_{\RR^2}|\widehat f(\xi)|^2(1+\vartheta(|\xi|))\,d\xi <\infty\}.
\end{equation}
\end{definition}
\noindent
We are interested in a $\vartheta$ that is associated with a modulus of continuity $\omega$ as follows:
\begin{definition}
For a given modulus of continuity $\omega(r)$, let us define
\begin{equation}\label{def:Omega}
\vartheta(r)=\begin{cases} \int_1^r\frac{ds}{s\,\om^2(s/r)} & \hbox{for}\ r>1 \\
0  & \hbox{for}\ 0\leq r\leq 1. \end{cases}
\end{equation}
\end{definition}
\noindent
For example, in the H\"older case $\om(r)=r^\gamma$ for $\gamma\in (0,1)$, we have $\vartheta(r)\approx c\,r^{2\gamma}$ as $r\to\infty$ and $W^{\vartheta,2}(\RR^2)$ coincides with the fractional-order Sobolev space $H^{\gamma,2}(\RR^2)$ defined as Bessel potentials. 
On the other hand, for $\varpi(r)$ as in \eqref{eq:bar-om}, we get
\begin{equation}\label{DiniExample2}
\vartheta(r)=\int_1^r \frac{ds}{s\,(\log(r/s))^{-2\alpha}}\approx \frac{|\log r|^{2\alpha+1}}{2\alpha+1}\ \hbox{as}\ r\to \infty.
\end{equation}
Note that $H^{\gamma,2}(\RR^2)\subset W^{\vartheta,2}(\RR^2)$ for any $\gamma\in (0,1)$.

Now let us explain why we are interested in $W^{\vartheta,2}(\RR^2)$.
In the H\"older case $\om(r)=r^\gamma$, functions in $H^{\gamma,2}(\RR^2)$ coincide with functions in $L^2(\RR^2)$ for which the $L^2$-modulus of continuity
\begin{equation}\label{def:L2-moduluscontinuity}
M_2(f,y)=\|f(\cdot + y)-f(\cdot)\|_{L^2}
\end{equation}
is small enough as $|y|\to 0$ that $M_2(f,y)\,|y|^{-1-\gamma}\in L^2(\RR^2)$; cf.\ \cite{S}.
Functions in $W^{\vartheta,2}$ can be similarly  characterized;
 for this purpose we need to introduce
\begin{equation}\label{def:tilde-omega}
\widetilde\omega(r)=\begin{cases} \omega(r) & \text{for $0<r<1$} \\ \frac{1}{\omega(1/r)} & \text{for $r>1$.}\end{cases}
\end{equation}
Note that $\widetilde\om(r)\to\infty$ as $r\to\infty$. In fact, for any $\alpha>0$ we have
\[
\omega(r)=|\log r|^{-\alpha}\quad \hbox{for $0<r\leq1/2$} \quad\Rightarrow\quad 
\widetilde\omega(r)=(\log r)^\alpha \quad \hbox{for $r\geq 2$.}
\]
However, in order to characterize functions in $W^{\vartheta,2}$ in terms of the modulus of continuity $\omega$, we
need to assume that $\omega$ satisfies the  ``square-Dini condition''
 \begin{equation}\label{def:squareDini}
 \int_{0}^\e \frac{(\omega(r))^2}{r}\,dr  < \infty \quad \hbox{for some } \e > 0.
 \end{equation}

\begin{lemma}\label{le:NASC-W^Om} 
If $\omega$ satisfies \eqref{def:squareDini}, then $f\in W^{\vartheta,2}$ if and only if $f\in L^2$ and 
\begin{equation}\label{eq:NASC-W^Om}
\int_{\RR^2}\frac{(M_2(f,y))^2}{|y|^2\,\widetilde \om^2(|y|)}\,dy<\infty.
\end{equation}
Here $\vartheta$ is defined in terms of $\omega$ by \eqref{def:Omega}.
\end{lemma}

\medskip\noindent
{\bf Proof.} Proceeding as in \cite{S}, we first use the Plancherel Theorem to obtain
\[
M_2^2(f,y)=(M_2(f,y))^2=\int_{\RR^2}|\widehat f(\xi)|^2|e^{-2\pi i\xi\cdot y}-1|^2\,d\xi.
\]
 Now let us consider the integral
\[
\int_{\RR^2}\frac{M_2^2(f,y)}{|y|^2\,\widetilde\omega^{2}(|y|)}\,dy=\int_{\RR^2}|\widehat f(\xi)|^2\,I(\xi)\,d\xi,
\]
where
\begin{equation}\label{def:I}
I(\xi)=\int_{\RR^2}\frac{|e^{-2\pi i\xi\cdot y}-1|^2}{|y|^2\,\widetilde\om^{2}(|y|)}\,dy.
\end{equation}
We note that the integral defining $I(\xi)$ converges: write 
\[
\int_{\RR^2}\frac{|e^{-2\pi i\xi\cdot y}-1|^2}{|y|^2\,\widetilde\om^{2}(|y|)}dy
=\int_{|y|<1}\frac{|e^{-2\pi i\xi\cdot y}-1|^2}{|y|^2\,\om^{2}(|y|)}dy+
\int_{|y|>1}\frac{|e^{-2\pi i\xi\cdot y}-1|^2\,\om^{2}(|y|^{-1})}{|y|^2}dy.
\]
Use $|e^{-2\pi i\xi\cdot y}-1|\leq c|y|$ for $|y|<1$ and $\om(r)\geq cr^\e$ for $0<r<1$ to conclude that the first integral  converges, and for the second integral use
 $|e^{-2\pi i\xi\cdot y}-1|\leq 2$  with
 \[
 \int_{|y|>1}\frac{\om^{2}(|y|^{-1})}{|y|^2}\,dy=c\,\int_1^\infty \frac{\om^{2}(r^{-1})}{r}\,dr=c\,\int_0^1 \frac{\om^{2}(s)}{s}\,ds<\infty.
 \]
Moreover, since $I(\xi)$ is rotation-invariant we can define the radial function $I_0(r)$:
\begin{equation}\label{def:I_0}
I_0(|\xi|):=I(\xi).
\end{equation}
The proof is complete if we can show that  $\vartheta(r)\approx I_0(r)$ as $r\to\infty$,
i.e.\ there exists a constant $c > 0$ such that
\begin{equation}\label{Claim}
c\,I_0(r)\leq \vartheta(r)\leq (1/c)\, I_0(r)\quad\hbox{for $r$ sufficiently large.}
\end{equation}
To estimate $I_0(r)$ as $r\to\infty$, let us choose $\xi=r\,(1,0)$ and let  $\tilde y=ry$. Then $I_0(r)=I_1(r)+I_2(r)$ where
\[
I_1(r)=\int_{|y|<1}\frac{|e^{-2\pi ir y_1}-1|^2}{|y|^2\,\om^{2}(|y|)}\,dy=
\int_{|\tilde y|<r} \frac{|e^{-2\pi i \tilde y_1}-1|^2}{|\tilde y|^2\,\om^{2}(|\tilde y|/r)}\,d\tilde y
\]
and
\[
I_2(r)=\int_{|y|>1}\frac{|e^{-2\pi i ry_1}-1|^2\,\om^{2}(|y|^{-1})}{|y|^2}\,dy=
\int_{|\tilde y|>r} \frac{|e^{-2\pi i  \tilde y_1}-1|^2\,\om^{2}(\frac{r}{|\tilde y|})}{|\tilde y|^2}\,d\tilde y.
\]
Recall that $\vartheta(r)\to\infty$ whereas $I_2(r)$ is decreasing as $r\to\infty$, so we need only concern ourselves with $I_1(r)$.
Moreover, it is clear that $I_1(r)-I_1(1)\leq c\,\vartheta(r)$ since
\[
\int_{1<|\tilde y|<r} \frac{|e^{-2\pi i \tilde y_1}-1|^2}{|\tilde y|^2\,\om^{2}(|\tilde y|/r)}\,d\tilde y\leq 4\int_1^r\frac{d\rho}{\rho\,\om^{2}(\rho/r)}=4\,\vartheta(r),
\]
so we need only show $I_1(r)\geq c\,\vartheta(r)$ for some $c>0$.  

To show $I_1(r)\geq c\,\vartheta(r)$, let us write $\tilde y_1=s\cos\theta$ and compute
$ |e^{-2\pi i \tilde y_1}-1|^2=2(1-\cos[2\pi s \cos\theta]).$
Consequently, 
\[
\begin{aligned}
I_1(r)=&\int_0^r  \int_{-\pi}^\pi \frac{2(1-\cos[2\pi s \cos\theta])}{s\,\om^2(s/r)} d\theta ds \\
&\geq 2\int_1^r \frac{1}{s\,\om^2(s/r)} \left(\int_{-\pi}^\pi (1-\cos[2\pi s \cos\theta])d\theta\right)ds.
\end{aligned}
\]
Thus it suffices to have
$
\int_{-\pi}^\pi (1-\cos[2\pi s \cos\theta]) d\theta \geq c \quad \hbox{for all}\ s\geq 1.
$
This is proved in Lemma \ref{technical-lemma} in the Appendix.
 $\Box$
 
 \medskip
 Since $\vartheta(r)$ is increasing in $r>1$, it is elementary to verify the following:
 \begin{lemma}\label{stom1} 
If $f\in W^{\vartheta,2}(\RR^2)$, $R_0>1$, and $0\leq\nu\leq 1$, then $\widehat f$ satisfies
\[
\int_{|\xi|\geq R_0}|\widehat f(\xi)|^2\,\vartheta(|\xi|)^\nu\,d\xi\leq \frac{\|f\|^2_{W^{\vartheta,2}}}{\vartheta(R_0)^{1-\nu}}.
\]
\end{lemma}

\medskip
Now, for the given $\alpha>1$ in $\varpi$, we want to consider a modulus of continuity $\om$ satisfying
\begin{equation}\label{om}
\om(r)=|\log r|^{-\beta}   \quad\hbox{for}\ 0<r\leq 1/2, \ \hbox{where $0<\beta<\alpha$.}
\end{equation}
 Note that $\om$ is {\it weaker} than $\varpi$: $\varpi(r)\leq\om(r)$ for $0<r\leq 1/2$ (and hence for all $r>0$).
We not only want $\om$ to satisfy the Dini condition, but we want $\varpi/\om$ to satisfy the square-Dini condition \eqref{def:squareDini} at $r=0$. In fact, if we assume that
 \begin{equation}\label{condition-a,b}
 1 < \beta < \alpha - 1/2,
 \end{equation}
which is possible since $\alpha>3/2$, then we have
 \begin{equation}\label{Condition-C_0}
 \mathcal{C}_{\alpha,\beta}:= \int_0^\infty \frac{\varpi^2(r)}{r\,\widetilde\om^2(r)}\,dr <\infty.
 \end{equation}
We will also use the notation $\mathcal{C}_{\alpha,\beta}$ for $K \,\mathcal{C}_{\alpha,\beta}$ where $K$ is a constant that might depend on other parameters. As in \eqref{DiniExample2} we find that for any constant $c>0$
\begin{equation}\label{est:vartheta}
\vartheta(r)\approx \vartheta(c\,r)\approx \frac{|\log r|^{2\beta+1}}{2\beta+1} \quad\hbox{as}\ r\to\infty.
\end{equation}
In particular, we have
\begin{equation}\label{LowerBound-vartheta}
 \vartheta(r)\geq C\,|\log r|^{3+\delta} \quad \hbox{for}\ r\geq 2,
 \end{equation}
 for $\delta=2(\beta-1)>0$ and some  $C>0$.

\begin{lemma}\label{stom} 
Suppose $\mu\in C_0^\varpi(\DD)$. 
\begin{enumerate}
\item[(a)]  For $1\leq p<\infty$, the $L^p$-modulus of $\mu$ is uniformly bounded by $\varpi$: 
\[
M_p(\mu,y):=\|\mu(\cdot +y)-\mu(\cdot)\|_{L^p}\leq 
C\,\varpi(|y|)\,\|\mu\|_{C^\varpi}\quad\hbox{for any $y\in\RR^2$}.
\]
\item[(b)] For $\om$ as in \eqref{om},\eqref{condition-a,b} with associated $\vartheta$, we have $\mu\in W^{\vartheta,2}(\RR^2)$ and
\[
\|\mu\|_{W^{\vartheta,2}}\leq \mathcal{C}_{\alpha,\beta}\,\|\mu\|_{C^\varpi}.
\]
\item[(c)] If $\om$ and $\vartheta$ are as in (b) and $f\in W^{\vartheta,2}(\RR^2)$, then $\mu\,f\in W^{\vartheta,2}(\RR^2)$ and
\[
\|\mu\,f\|_{W^{\vartheta,2}}\leq  \mathcal{C}_{\alpha,\beta}\,\|\mu\|_{C^\varpi}\,\|f\|_{W^{\vartheta,2}}.
\]
\end{enumerate}
\end{lemma}

\medskip\noindent
{\bf Proof.} To prove (a), we begin with $|\mu(x+y)-\mu(x)|\leq\varpi(|y|)\,\|\mu\|_{C^\varpi}$.
For $|y|\leq 1$,
\[
\begin{aligned}
M^p_p(\mu,y)=\int_{\DD_2}|\mu(x+y)-\mu(x)|^p\,dx 
&\leq \varpi^p(|y|)\,\|\mu\|^p_{C^\varpi}\int_{\DD_2}\,dx.
\end{aligned}
\]
For $|y|\geq 1$ we can use $\varpi(|y|)=\varpi(1/2)$ to conclude
\[
M_p(\mu,y)\leq 2\|\mu\|_{L^p}\leq 2\|\mu\|_{\infty}\left(\int_{\DD_1}\,dx\right)^{1/p}
\leq C_{p}\,\varpi(|y|)\, \|\mu\|_{C^\varpi},
\]
where $\|\cdot\|_\infty$ denotes the sup-norm.

To prove (b), we use (a) and some of the formulas in the proof of Lemma \ref{le:NASC-W^Om}.
We first use (\ref{Claim}) and then (a) to estimate
\[
\begin{aligned}
\int_{|\xi|>1} |\widehat\mu(\xi)|^2\vartheta(|\xi|)\,d\xi & \leq C\int_{\RR^2} |\widehat\mu(\xi)|^2 I(\xi)\,d\xi 
 =C \int_{\RR^2} \frac{M_2^2(\mu,y)}{|y|^2\widetilde\om^{2}(|y|)}\,dy \\
&\leq C\, \int_{\RR^2} \frac{\varpi^2(|y|)}{|y|^2\widetilde\om^{2}(|y|)}dy\, \|\mu\|^2_{C^\varpi}
= \mathcal{C}_{\alpha,\beta}\|\mu\|^2_{C^\varpi}.
\end{aligned}
\]

\smallskip
To prove (c), note that
\begin{equation}\label{Bel1}
\left\|\mu\,f\right\|^2_{W^{\vartheta,2}}\approx \int_{\RR^2}\frac{M_2^2(\mu f,y)}{|y|^2\,\widetilde\om^{2}(|y|)}\,dx+ \left\|\mu\,f\right\|^{2}_{L^2}
\end{equation}
But we have
\[
\begin{aligned}
M_2^2(\mu f,y)&=\int_{\RR^2}{|(\mu(z+y)-\mu(z))\,f(z+y)+\mu(z)\,(f(z+y)-f(z)))|}^2 \,dz\\
&\leq C\left({\varpi^2(|y|)\,\left\|\mu\right\|^2_{C^\varpi}\,\left\|f\right\|^{2}_{L^2}+\left\|\mu\right\|^2_{\infty}\,M_2^2\,(f,y)}\right).
\end{aligned}
\]
Hence (\ref{Bel1}) gives,\\
\[
\begin{aligned}
\left\|\mu\,f\right\|^2_{W^{\vartheta,2}} 
 & \leq C{ \int_{\RR^2}\frac{\varpi^2(|y|)\,\left\|\mu\right\|^2_{C^\varpi}\,\left\|f\right\|_{L^2}^2+\left\|\mu\right\|^{2}_{L^\infty}\,M_2^2\,(f,y)}{|y|^2\,\widetilde\om^{2}(|y|)}\,dy+ \left\|\mu\right\|_{{\infty}}^2\,\left\| f \right\|^2_{L^2}}\\
& \leq C\left\|\mu\right\|_{C^{\varpi}}^2\,{ \int_{\RR^2}\frac{\varpi^2(|y|)\,\left\|f\right\|_{L^2}^2+ M_2^2\,(f,y)}{|y|^2\,\widetilde\om^{2}(|y|)}\,dy+ \left\|\mu\right\|^2_{C^{\varpi}}\,\left\|f\right\|^2_{L^2}}\\
&\leq \mathcal{C}_{\alpha,\beta}\,\left\|\mu\right\|_{C^{\varpi}}^2\,\left\|f\right\|_{W^{\vartheta,2}}^2.
\end{aligned}
\]
This completes the proof. \ $\Box$


\section{Estimating the complex geometric optics solutions as $|k|\to \infty$} \label{Estimates-k}

The proof of the Main Theorem uses the so-called ``complex geometric optics solutions'' which were initiated by Calderon \cite{C}.\ As in \cite{AP}, these can be constructed by solving the associated  $\RR$-linear Beltrami equation
\begin{equation}\label{def:BeltramiEqn}
\overline{\del}f=\mu\, \overline{\del f}, \quad\hbox{where}\ \mu=(1-\gamma)/(1+\gamma).
\end{equation}
Here $\overline{\del}=\del_{\overline z}=\frac{1}{2}(\del_x+i\del_y)$ and $\del=\del_z=\frac{1}{2}(\del_x-i\del_y).$
In fact,  a complex solution $u_\gamma$ of  (\ref{eq:conductivity}) can be recovered from a solution $f_\mu$ of (\ref{def:BeltramiEqn}) simply by letting
\begin{equation}\label{eq:f->u}
u_\gamma=\operatorname{Re}(f_\mu)+i\, \operatorname{Im}(f_{-\mu}).
\end{equation}
As in \cite{AP} and \cite{BFR} (see also Theorem \ref{WhtineyExt} below), the problem can be reduced to the case that $U=\DD$ and $\gamma=1 $ in a neighborhood of $\partial\DD$, so $\mu$ has compact support in $\DD$ and we may consider \eqref{def:BeltramiEqn} on $\CC$. Henceforth we use the notation $\|\cdot\|_\infty$ for the sup-norm on $\CC$ (or $\RR^2$).

Now, for each $k\in\CC$, \cite{AP} shows that there is a unique solution $f_\mu(z,k)$ of (\ref{def:BeltramiEqn}) in the form
\begin{equation}\label{geom-optics-condn1}
f_\mu(z,k)=e^{ik\phi(z,k)},
\end{equation}
where, for each fixed $k\in\CC$, $\phi(z,k)$ is a quasiconformal\footnote{A homeomorphism $\phi:\CC\to\CC$ is $K$-{\it quasiconformal} if it is orientation-preserving, it $\phi\in H^{1,2}_{\loc}(\CC)$, and the directional derivatives are $\partial_\nu\phi$ satisfy the $\max_{\nu}|\partial_\nu\phi(z)|\leq K\min_\nu|\partial_\nu\phi(z)|$ for almost every $z\in\CC$. If $\phi$ is $K$-quasiconformal, then it is locally $K^{-1}$-Holder continuous. See  \cite{AIM}.} homeomorphism in $z$ and satisfies the nonlinear Beltrami equation
\begin{equation}\label{nonlinearBeltrami}
\del_{\overline{z}}\phi=-\frac{\overline k}{k}\mu(z)e_{-k}(\phi(z))\,\overline{\del_z\phi},
\end{equation}
where $e_{k}(z)=\exp(i(kz+\overline{k}\overline{z}))$, and the boundary condition
\begin{equation}\label{BC-phi}
\phi(z)=z+O(z^{-1})\quad\hbox{as}\ |z|\to\infty.
\end{equation}
 In \eqref{nonlinearBeltrami} and \eqref{BC-phi} we consider $k$ fixed, but the dependence of $\phi$ on $k$ is important. For $\mu\in L^\infty$, it was shown in \cite{AP} that $\phi(z,k)\to z$ as $|k|\to \infty$ uniformly for $z\in \CC$: 
 \begin{equation}\label{phi-large-k-2}
|\phi(z,k)- z| \leq C(k),
\end{equation}
where $C(k)\to 0$ as $|k|\to\infty$. For more regular conductivities,  $C(k)$ can be described more precisely. For H\"older continuous conductivities $\mu\in C^\gamma$, it was shown in \cite{BFR} that \eqref{phi-large-k-2} holds with $C(k)=c\,|k|^{-a}$, where $c,a>0$ depend on $\gamma$ (and other parameters). For the Dini continuous conductivities that we consider, $\mu\in C^\varpi$,  $C(k)$ will depend on $\om$ and $\vartheta$ as in the previous section.

As in \cite{AP}  and \cite{BFR}, the study of the solutions  of (\ref{nonlinearBeltrami})-(\ref{phi-large-k-2}) is reduced to  the study of solutions $\psi(z,k)$ of a $\CC$-linear Beltrami equation with boundary condition:
\begin{equation}\label{linearBeltrami}
\begin{aligned}
\del_{\overline{z}}\psi&=-\frac{\overline k}{k}\mu(z)e_{-k}(z)\,{\del_z\psi},\\
\psi(z,k)&=z+O(z^{-1})\quad\hbox{as}\ |z|\to\infty.
\end{aligned}
\end{equation}
By writing $\psi=z+\eta$, this is reduced to solving a nonhomgenous Beltrami equation with boundary condition
$\eta=O(z^{-1})$ as $|z|\to\infty$. The function $\eta$ is uniquely determined in $H^{1,p}(\CC)$ for all $p>2$,
so  $\psi\in H^{1,2}_{\loc}(\CC)$  is uniquely determined.

A tool used in the study of \eqref{linearBeltrami} is the Beurling transform
\begin{equation}
T(g)(z)=-\frac{1}{\pi}\int \frac{g(w)}{(w-z)^2}\,dw,
\end{equation}
which has the well-known properties (cf.\ \cite{AIM}) that 
$T:L^p(\CC)\to L^p(\CC)$ is bounded for all $1<p<\infty$ with
\begin{equation}\label{Lp-normofT}
\lim_{p\to 2}\|T\|_{L^p\to L^p}=\|T\|_{L^2\to L^2}=1,
\end{equation}
and $T(\del_{\overline{z}}\eta)=\del_z \eta$ for all $\eta\in H^{1,p}(\CC)$.
As in \cite{AP}, there is a convergent series representation
\begin{equation}\label{seriesrepresentation1}
\del_{\overline{z}}\psi=\sum_{n=0}^\infty (aT)^n a, \quad\hbox{where} \quad a(z,k)=-\frac{\overline k}{k}\mu(z)e_{-k}(z),
\end{equation}
provided $\|aT\|_{L^p\to L^p}<1$. Since $|a(z,k)|\leq \|\mu\|_\infty\leq \kappa<1$,   this condition  
will be met provided $p$ is sufficiently close to $2$. Now choose $p=p(\kappa)>2$ so that
\begin{equation}\label{p-condition}
 \kappa_1:=\kappa\|T\|_{L^p\to L^p}  <1\quad\hbox{for all $p\in[2,p(\kappa)]$}.
\end{equation}
Then (\ref{seriesrepresentation1}) converges in $L^p$ for all $p\in[2,p(\kappa)]$.
For a chosen integer $n_0>0$ , write
\begin{equation}\label{dopsidef}
\del_{\overline{z}}\psi=g_k(z)+h_k(z):= \sum_{n=0}^{n_0-1} (a(z,k)T)^n\,a(z,k)+\sum_{n=n_0}^\infty (a(z,k)T)^n\,a(z,k).
\end{equation}
For $2\leq p\leq p(\kappa)$, a geometric series can be used to estimate
\begin{equation}\label{h-est}
\sup_{k\in\CC}\|h_k\|_{L^p}\leq \pi^{1/p}\,\kappa\,\frac{\kappa_1^{n_0}}{1-\kappa_1}.
\end{equation}
By choosing $n_0$ sufficiently large,  $\sup_k \|h(\cdot,k)\|_{L^p}$ can be made arbitrarily small.
Similarly, we can estimate
\begin{equation}\label{g-est}
\sup_{k\in\CC}\|g_k\|_{L^p}\leq\pi^{1/p}\,\kappa\,\frac{1-\kappa_1^{n_0}}{1-\kappa_1}\leq C_1(\kappa).
\end{equation}
We cannot make this small so we will require estimates on the Fourier transfrom of $g$.
The following result is analogous to Lemma 3.6(c) in \cite{BFR}:

\begin{lemma}\label{decomp-g+h} 
Suppose $\mu\in C_0^\varpi(\DD)$ with  $\|\mu\|_\infty\leq\kappa$ and $\|\mu\|_{C^\varpi}\leq \Gamma$. Let $\om$ satisfy \eqref{om},\eqref{condition-a,b} with associated $\vartheta$ as in \eqref{def:Omega}.
Let $\psi\in H^{1,2}_{\loc}(\CC)$ be the unique solution of (\ref{linearBeltrami}) and  for any fixed $n_0\in\NN$ consider the decomposition \eqref{dopsidef}.
Then, for any $R_0>1$ and $|k|\geq 2R_0$, we can estimate
\begin{equation}\label{est-g-hat}
\int_{|\xi|<R_0}|\widehat g_k(\xi)|^2\,d\xi \leq \frac{n_0\,(\mathcal{C}_{\alpha,\beta}\Gamma)^{n_0}}{{\vartheta(|k|/2)}}.
\end{equation}
\end{lemma}

\medskip\noindent
{\bf Proof:} Write $g_k(z)$ as
\begin{subequations}
\begin{equation}\label{defg}
g_k(z)=\sum_{n=0}^{n_0-1}(-\overline{k}\,/k)^{n+1}\,e_{-(n+1)\,k}(z)\,f_n(z),
\end{equation}
where $f_0 =\mu$ and
\begin{equation}\label{deff_n}
 f_n(z)=\mu\,T_n\,\mu\,T_{(n-1)}\,\mu.........\mu\,T_1(\mu)\quad\hbox{for } n>0. 
 \end{equation}
\end{subequations}
Here $ T_j =e_{jk} T e_{-jk}$ is the Fourier multiplier with  symbol $(\xi-jk)/\overline{(\xi-jk)}$, which is unimodular, so $ \|T_j\|_{W^{\vartheta,2}\rightarrow W^{\vartheta,2}}=1$ and 
$\|f_n\|_{W^{\vartheta,2}}\leq \|\mu\|^{n+1}_{W^{\vartheta,2}}\leq \mathcal{C}_{\alpha,\beta}^{n+1} \|\mu\|^{n+1}_{C^\omega}$ where we have also used Lemma \ref{stom}(c). But this means
\begin{equation}\label{est-g_k}
\|g_k\|_{W^{\vartheta,2}}\leq \sum_{n=0}^{n_0-1} \|f_n\|_{W^{\vartheta,2}}\leq n_0\,\mathcal{C}_{\alpha,\beta}^{n_0} \|\mu\|^{n_0}_{C^\omega}.
\end{equation}

Now let us turn to \eqref{est-g-hat}:
\[
\int_{|\xi|<R_0}|\widehat g_k(\xi)|^2\,d\xi \leq \sum_{n=0}^{n_0-1}\int_{|\xi|<R_0}|\widehat f_n(\xi -(n+1)\,k)|^2 \,d\xi
\]
and by Lemma \ref{stom1}
 \[
\int_{|\xi|<R_0}|\widehat f_n(\xi -(n+1)\,k)|^2 \,d\xi \leq  \int_{|\xi|>(n+1)|k|-R_0}|\widehat f_n(\xi)|^2\,d\xi
\leq \frac{\|f_n\|^2_{W^{\vartheta,2}}}{\vartheta((n+1)\,|k|-R_0)}.
\]
For  $|k|\geq 2R_0$, since $\vartheta(|x|)$ is an increasing function,
$\vartheta((n+1)\,|k|-R_0)\geq\vartheta(|k|/2)$,
so we obtain
\[
\int_{|\xi|<R_0}|\widehat g_k(\xi)|^2\,d\xi \leq
\sum_{n=0}^{n_0-1} \frac{(\mathcal{C}_{\alpha,\beta}\Gamma)^{n+1}}{\vartheta(|k|/2)}
\leq \frac{n_0\,(\mathcal{C}_{\alpha,\beta}\Gamma)^{n_0}}{{\vartheta(|k|/2)}}   \quad \Box
\]                                                                                         

\medskip
Now we want to use these estimates to control the behavior of $\psi(z,k)-z$ as $|k|\to\infty$. For this we need to use the Cauchy transform
\begin{equation}\label{def:P}
P[g](z)=-\frac{1}{\pi}\int_{\CC}\frac{g(w)}{w-z}\,dw,
\end{equation}
which has the well-known property that it acts as an inverse for $\del_{\overline z}$:
$P\del_{\overline z}\rho=\del_{\overline z}P\rho=\rho$
for $\rho\in C_0^\infty(\CC)$.
It is well-known (cf.\ \cite{AIM}) that for $2<p<\infty$, the Cauchy transform is bounded
\begin{equation}\label{PLp->C^alpha}
P:L^p(\CC)\to C^\alpha(\CC)\quad\hbox{where $\alpha=1-2/p$.}
\end{equation}
\begin{proposition}\label{pr:estimateforlinearBeltrami} 
Suppose $\mu\in C_0^\varpi(\DD)$ with $\|\mu\|_\infty\leq\kappa$ and $\|\mu\|_{C^\varpi}\leq \Gamma$.
Let $\omega$ and $\vartheta$ be as in Lemma \ref{decomp-g+h} and
let $\psi\in H^{1,2}_{\loc}(\CC)$ be the unique solution of (\ref{linearBeltrami}). Then there exist constants $C=C(\alpha,\beta,\kappa,\Gamma)$ and $a=a(\kappa,\Gamma)>0$ such that for all $z\in\CC$ we have
\begin{equation}\label{est:psi-z}
|\psi(z,k)-z|\leq \frac{C}{\vartheta(|k|)^{a}}\quad\hbox{as $|k|\to\infty$.}
\end{equation}
\end{proposition}

\medskip\noindent
{\bf Proof.} We observe that $\del_{\overline z}\psi$ has compact support in $\DD$, so we may use \eqref{dopsidef} to write for fixed $k\in\CC$:
\[
\psi(z,k)-z=P[\del_{\overline z}\psi](z,k)
=P[g_k](z)+P[h_k](z).
\]
Let us choose $p=p(\kappa)$ as in \eqref{p-condition} and let $q=p'$ be the conjugate index; in particular, we have $p>2$ so $1<q<2$. Let us take
\begin{equation}\label{epsilon}
\e=\frac{C}{\vartheta(R_0)^a}\approx \frac{C}{(\log R_0)^{(2\beta+1)a}},
\end{equation}
where $R_0, a>0$ are to be determined so that $\|P[h_k]\|_\infty <\e/3$ and 
$\|P[g_k]\|_\infty <2\e/3$.

\smallskip\noindent{\bf Estimate $P[h]$.} Using \eqref{h-est}, we have
\[
\|P[h_k]\|_\infty \leq C(\kappa,q)(\kappa_1)^{n_0},
\]
so we will have $\|P[h_k]\|_\infty<\e/3$ provided
\begin{equation}\label{n_0}
n_0\geq \frac{\log(C_0\e)}{\log(\kappa_1)}\quad\hbox{where}\ C_0=C_0(\kappa,p).
\end{equation}

\smallskip\noindent{\bf Estimate $P[g]$.} Since $g_k(z)$ has compact support in $\DD$,  we can write
\[
P[g_k](z)=\int_{\CC}K_z(y)g_k(y)\,dy \quad\hbox{where}\ K_z(y)=\frac{1}{\pi}\frac{\chi(|y|)}{z-y};
\]
here $\chi(r)$ is a smooth cut-off function satisfying $\chi(r)=1$ for $0\leq r\leq 1$ and $\chi(r)=0$ for $r\geq 3/2$.
Notice that $K_z\in L^q(\CC)$ for all $1\leq q<2$. We want to use the Fourier transform to 
estimate $P[g_k]$.  The Fourier transform of $g_k$ is well-behaved since $g_k\in L^p$ and has compact support, so $g_k\in L^2$ and hence $\widehat g_k\in L^2$.
We also know that $\widehat{K_z}$ is in $L^p(\CC)$ for all $p>2$ (by the Hausdorff-Young inequality) so $\widehat{K_z}\in L^2(\DD)$; but we do {\it not} know that $\widehat{K_z}\in L^2(\CC)$, so we cannot just use Plancherel's theorem to conclude
\begin{equation}\label{eq:Plancherel}
P[g_k](z)=\int_{\CC} \widehat{K_z}(\xi)\,\widehat{g_k}(\xi)\,d\xi.
\end{equation}
However, the integral in \eqref{eq:Plancherel} converges by the H\"older inequality (since $\widehat K_z\in L^p(\CC)$ and $\widehat g_k\in L^q(\CC)$); then we can use an approximation argument to show its equality with $P[g_k](z)$. Now, for $R_0>0$ let us write
\[
\left| P[g_k]\right(z)| \leq 
\int_{|\xi|>R_0}  \left| \widehat{K_z}(\xi)\,\widehat{g_k}(\xi)\right|\,d\xi + \int_{|\xi|<R_0} \left| \widehat{K_z}(\xi)\,\widehat{g_k}(\xi)\right|\,d\xi=:I_1+I_2.
\]
Want $R_0$ large so that $I_1<\e/3$ and $I_2<\e/3$ for all $|k|\geq R_0$ 
  (uniformly in $z$).

To estimate $I_1$, we use the following estimate that is proved in the Appendix:
\begin{equation}\label{est:hat-K}
\left| \widehat{K_z}(\xi)\right|\leq C\,\frac{\log |\xi|}{|\xi|} \quad\hbox{for}\ |\xi|>2\quad \hbox{(uniformly in $z\in\CC$).}
\end{equation}
Using \eqref{est:hat-K} and \eqref{LowerBound-vartheta}, we conclude that
$
\widehat K_z(\xi)/\vartheta^{1/2}(|\xi|)\in L^2(|\xi|>2)
$
since 
\[
\int_{|\xi|>2}\frac{ |\widehat K_z(\xi) |^2}{\vartheta(|\xi|)}\,d\xi \leq \int_{2}^{\infty}\frac{( \log r)^{2}}{r\,\vartheta(r)}\,dr \leq \int_{2}^{\infty}\frac{1}{r\, {(\log r)}^{1+\delta}}\,dr < \infty.
\]
\noindent
Now by Cauchy-Schwartz, \eqref{LowerBound-vartheta}, and \eqref{est-g_k}:
\[
\begin{aligned}
(I_1)^2 \leq 
 \|g_k\|^2_{W^{\vartheta,2}}\, \int_{|\xi|>R_0}\frac{ |\widehat K_z(\xi) |^2}{\vartheta(|\xi|)}\,d\xi 
 \leq n_0\,\left({\mathcal{C}_{\alpha,\beta}}\, \Gamma\right)^{n_0}
 (\log R_0)^{-\delta},
  \end{aligned}
\]
for any $\delta\leq 2(\beta-1)$.
So we can choose $R_0$ large enough that $I_1<\e/3$ for all $|k|\geq R_0$.


To estimate $I_2$, we use Cauchy-Schwartz, \eqref{est:hat-K}, \eqref{est:vartheta}, and Lemma \ref{decomp-g+h}:
\[
\begin{aligned}
(I_2)^2 & \leq \int_{|\xi|<R_0} |\widehat K_z(\xi)|^2\,d\xi\, \int_{|\xi|<R_0} |\widehat g_k(\xi)|^2\,d\xi \\
& \leq C\left(\int_{|\xi|<1} |\widehat K_z(\xi)|^2\,d\xi + \int_1^{R_0}\frac{(\log r)^2}{r}\,dr\right)
\left(\frac{n_0\,(\mathcal{C}_{\alpha,\beta}\Gamma)^{n_0}}{{\vartheta(|k|)}}\right) \\
& \leq C\left(1+(\log R_0)^3\right) \frac{(\mathcal{C}_{\alpha,\beta}\Gamma)^{n_0+1}}{{\vartheta(|k|)}}.
\end{aligned}
\]
So, by \eqref{LowerBound-vartheta} we can choose $R_0$ sufficiently large that $I_2<\e/3$ for all $|k|\geq R_0$. 

Finally, we need to confirm that the choices of $n_0$ and $R_0$ are compatible. 
First, let $A:=2 {\mathcal C}_{\alpha,\beta}\Gamma$ so that $n_0 ({\mathcal C}_{\alpha,\beta}\Gamma)^{n_0}\leq A^{n_0}$ and then take equality in \eqref{n_0}. We find 
\[
n_0({\mathcal C}_{\alpha,\beta}\Gamma)^{n_0}\leq B\,(\log R_0))^{\tau(2\beta+1)a},\quad\hbox{where}\ \tau=\frac{\log A}{\log \kappa_1^{-1}}.
\]
(Here $B$ is also independent of $n_0,R_0,\e$.)
If we use this in our estimate for $I_1$ we find
$
I_1^2\leq C(\log R_0)^{\tau(2\beta+1)a-\delta},
$
where $\delta\leq 2(\beta-1)$.
So we can achieve $I_1<\e/3=C(\log R_0)^{-(2\beta+1)a}$ provided
\begin{equation}\label{cond:a}
a<\frac{\beta-1}{(2\beta+1)(1+\tau/2)}
\end{equation}
Similarly, we have $I_2^2\leq C\,(\log R_0)^{3+(2\beta+1)(\tau a-1)}$. We find that  \eqref{cond:a}  is exactly the condition we need to make $I_2<\e/3$.  This completes the proof. \ 
 $\Box$

\bigskip
We now want to obtain an estimate like \eqref{est:psi-z} for solutions $\phi(z,k)$ of the nonlinear Beltrami equation (\ref{nonlinearBeltrami})-(\ref{BC-phi}). This analysis uses the fact that for each $k\in\CC$, $\phi(z,k)$ is a quasiconformal homeomorphism, and hence has an inverse function $\psi:\CC\to\CC$ defined by
\begin{equation}\label{def:inverse-of-phi}
\psi\circ\phi(z)=z
\end{equation}
which is also quasiconformal. Differentiating (\ref{def:inverse-of-phi}) with respect to $z$ and $\overline z$ shows that $\psi$ satisfies
\[
\begin{aligned}
\del_{\overline{z}}\psi&=-\frac{\overline k}{k}\mu(\psi(z,k) )e_{-k}(z)\,{\del_z\psi},\\
\psi(z,k)&=z+O(z^{-1})\quad\hbox{as}\ |z|\to\infty.
\end{aligned}
\]
Of course, this is of the form (\ref{linearBeltrami}). So, provided we can show that the coefficient $\mu(\psi(z,k) )$ satisfies the conditions of Proposition \ref{pr:estimateforlinearBeltrami}, we may use it to conclude the desired estimate for $\psi$; these estimates then apply to $\phi=\psi^{-1}$. We use this line of reasoning to prove the following:

\begin{theorem}\label{th:estimatefornonlinearBeltrami} 
Suppose $\mu\in C_0^\varpi(\DD)$ and $\|\mu\|_{C^\varpi}\leq \Gamma$.
Let $\omega$ and $\vartheta$ be as in Lemma \ref{decomp-g+h} and
let  $\phi\in H^{1,2}_{\loc}(\CC)$ be the unique solution of (\ref{nonlinearBeltrami})-(\ref{BC-phi}). Then there exist positive constants $C_*=C_*(\kappa,\Gamma)$ and $a=a(\kappa,\Gamma)$ such that
\[
|\phi(z,k)-z|\leq \frac{C_*}{\vartheta(|k|)^{a}}\quad\hbox{as $|k|\to\infty$.}
\]

\end{theorem}
\medskip\noindent
{\bf Proof.} Let $\widetilde\mu(z,k)=\mu(\psi(z,k) )$. As in \cite{AP}, since $\psi$ is H\"older continuous and $\mu$ has support in $\DD$, $\widetilde\mu$ has support in $\DD_4$ (by the 1/4-Koebe theorem) for each fixed $k\in\CC$. So we can apply Proposition \ref{pr:estimateforlinearBeltrami} in $\DD_4$ provided we can show $\widetilde\mu\in C^{\varpi}_0(\DD_4)$. Let $\gamma$ be the H\"older coefficient for $\psi$. (Recall that $\gamma=K^{-1}$ where $K=(1-\kappa)/(1+\kappa)$.) Then
\[
\frac{|\widetilde\mu(z)-\widetilde\mu(y)|}{\varpi(|z-y|)}\leq \|\mu\|_{C^\varpi} \sup_{z,y}\frac{\varpi(|\psi(z)-\psi(y)|)}{\varpi(|z-y|)}.
\]
But we know $|\psi(z)-\psi(y)|\leq C_1\,|z-y|^\gamma$ for all $z,y\in \DD_4$, where we may assume $C_1>1$.
Since $\varpi(r)$ is nondecreasing, we have $\varpi(|\psi(z)-\psi(y)|)\leq \varpi(C_1|z-y|^\gamma)$, so it suffices to show that
\begin{equation}\label{tilde-mu-condition}
\sup_{0<r<\infty}\frac{\varpi(C_1\,r^\gamma)}{\varpi(r)}<\infty.
\end{equation}
But since we know explicitly that $\varpi(r)=|\log r|^{-\alpha}$ near $r=0$ and is constant for large $r$, condition \eqref{tilde-mu-condition} is easily verified.
$\Box$


\section{Regularity of the complex geometric optics solutions }\label{sec:regularity}

Now let us turn to the regularity of the complex geometric optics solutions $f_\mu$ and $u_\gamma$. Similar to the result obtained
in \cite{BFR}, we show that the assumption $\|\mu\|_{C^\varpi}\leq \Gamma$ gives 
an upper bound for $f_\mu$ and lower bound for the Jacobian, $J_{f_\mu}$. We obtain interior estimates on compact subsets which can then be used to prove the main theorem when $\mu_j $ has compact support in $\DD$.
To obtain these estimates, we first obtain estimates for a more general Beltrami equation which is also $\RR$-linear. To analyse such an equation, we use a modulus of continuity $\si$ as  defined in \cite{MT}.

\medskip\noindent
{\bf Definition of $\si(r)$:}
 Define a continuous function $\si$ on $[0,\infty)$ by
\begin{equation}\label{defsigma}
 \si(r) :=\int_0^{r} \frac{\varpi(s)}{s}\,ds= \frac{|\log r|^{1-\alpha}}{\alpha-1}
 \quad\hbox{for}\ 0< r\leq 1/2,
\end{equation}
with $\si(0)=0$ and $\si(r)$ is a constant for $r\geq 1/2$.
(Note that $\sigma$ need not satisfy the Dini condition at $r=0$.)


\medskip\noindent 
Using $\si$ as the modulus of continuity, 
$C^\si(\overline U)$ and $C^{1,\si}(\overline U)$ are Banach spaces defined as in section 2. It is easy to see that $\varpi(r)\leq \si(r)$ and hence for any domain $\Om$, 
\begin{equation}\label{omsirel}
C^\varpi(\Om) \subset C^\si(\Om). 
\end{equation}

\begin{proposition}\label{bel}   
For a bounded domain $\Om$, let $\mu,\nu \in C^{\varpi}(\Om)$ satisfying $|\mu(\xi)| +|\nu(\xi)| \leq \kappa<1$ for all $\xi\in \Om$. Let $v\in H^{1,2}_{\loc}(\Om)$  be a solution to the equation
\begin{equation}\label{bequation}   
\overline{\del}v- \mu\,\del v- \nu\, \overline{\del v}=0.
\end{equation}
Consider domains $D$, $U$  such that $\overline D \subset U$ and $\overline U \subset \Om$. If $\|\mu\|_{C^{\varpi}(\overline U)} +\|\nu\|_{C^{\varpi}(\overline U)}<\Gamma$, then we have the following:
\begin{enumerate}
\item[(a)] 
 $v\in C^{1,\si}(\Om)$. In particular, we have $v\in C^{1,\si}(\overline D)$ and there exists $ K_1=K_1(\kappa,\Gamma,D, U)$ such that
\begin{equation}\label{finaltresult}
\| v\|_{C^{1,\si}(\overline D)}\leq K_1\, \|v\|_{C^{0}(\overline U)}.
\end{equation}
\item[(b)] 
If $v$ is a quasiconformal homeomorphism in $\CC$, let $M=M(U)$ satisfy
\begin{equation}\label{def:M}
M=\max_{x\in U} |v(x)|.
\end{equation}
Then there exists a constant $ K_2=K_2(\kappa,\Gamma,D, U,M)>0$ such that
\begin{equation}\label{jacobian}
 \inf_{z\in D} J_v(z)=\inf_{z\in D}(|\del v(z)|^2 - |\overline \del v(z)|^2)\geq  K_2.
\end{equation}
\end{enumerate}
\end{proposition}
\noindent 
{\bf Proof:} For the proof of (a), refer to \cite{BV}, where a similar result is proved for the non-homogeneous equation corresponding to (\ref{bequation}).

To prove (b), we proceed as in \cite{BFR}. For $z\in D$, we use \eqref{bequation} to obtain
\begin{equation}\label{jacobian-1}
 J_v(z)=|\del v(z)|^2 - |\overline \del v(z)|^2 \geq (1- \kappa^2)|\del v(z)|^2.  
\end{equation}
Consider the inverse function $v^{-1}$ of $\xi=v(z)$, which satisfies the Beltrami equation 
\[
\partial_\xi(v^{-1})-(\mu\circ v^{-1})\overline{\partial_\xi(v^{-1})}-(\nu\circ v^{-1})\partial_\xi (v^{-1})=0.
\]
  Since $v^{-1}$ is quasiconformal, it is H\"older continuous, so the coefficients $\mu \circ v^{-1}$ and $\nu \circ v^{-1}$ in this Beltrami equation are in ${C^\varpi(v(D))}$, by a similar argument as in the proof of Theorem \ref{th:estimatefornonlinearBeltrami}.\ Therefore $v^{-1}$ satisfies the conditions in (a) to obtain the corresponding estimate (\ref{finaltresult}). In particular, we have
\begin{equation}\label{vinvresult}
|\del_{\xi} v^{-1}\circ v(z)|\leq K_3 \quad\hbox{for $z\in D$,}
\end{equation}
where $K_3= K_3(\kappa,\Gamma,D, U,M)$.
On the other hand, differentiating $z=v^{-1}\circ v(z)$ by the chain rule, we obtain
$
1=|(\partial_\xi v^{-1}\circ v)(\partial_z v)|=|(\partial_\xi v^{-1}\circ v)|\,|(\partial_z v)|.
$
So, by \eqref{vinvresult}, we have
\begin{equation}\label{est:Dv}
|\partial_z v(z)|=\frac{1}{|\partial_\xi v^{-1}\circ v|}\geq \frac{1}{K_3}\quad\hbox{for $z\in D$},
\end{equation}
Since $J_v(z)\geq |\partial_z v(z)|^2$, we obtain \eqref{jacobian}.
  \hfill$\Box$\\

\noindent
We can now obtain the following result: 

\begin{theorem}\label{th:f-regularity}. 
Suppose $\mu\in C_0^\varpi(\DD)$ with  $\|\mu\|_{C^\varpi}\leq \Gamma$. There exist positive constants $C_1(\kappa,\Gamma,|k|)$ and $C_2(\kappa,\Gamma,|k|)$ so that the complex geometric optics solution (\ref{geom-optics-condn1}) satisfies
\begin{equation}\label{eq:regularity-f}
\|f_\mu(\cdot,k)\|_{C^{1,\sigma}(\DD)}\leq C_1 \quad\hbox{and}\quad \inf_{z\in\DD}|J_{f_\mu}(z,k)|\geq C_2.
\end{equation}
\end{theorem}
\noindent
 {\bf Proof:} Recall from (\ref{geom-optics-condn1}) that $ f_\mu(z,k)=e^{ik\phi(z,k)}$ where  $\phi = z+ \e(z,k)$, with $\e(z,k) $ uniformly bounded for fixed $k$,  is H\"older continuous. So the coefficient $\frac{\overline k}{k}\mu e_{-k}(\phi)$ in (\ref{nonlinearBeltrami}) is in $C^\varpi(\DD_3)$. Also, for $z\in \DD$, max $|\phi(z,k)|= C$ for $C=C(k, \kappa, \DD)$ gives us bounds for $ f_\mu$ as $1/C \leq |f_\mu (z,k)|\leq C$.
Hence we can apply Proposition \ref{bel} to $\phi(z,k)$ to obtain
$\| \phi\|_{C^{1,\si}(\overline \DD)}\leq K_1\, \|\phi\|_{C^{0}(\overline \DD_2)}$
 for some constant $K_1=K_1(\kappa,\Gamma)$.  This in turn shows that there exists a constant $C_1(\kappa,\Gamma,|k|)$ such that $\|f_\mu(\cdot,k)\|_{C^{1,\sigma}(\DD)}\leq C_1 $.

Proposition \ref{bel} can also be used to obtain the lower estimate 
\begin{equation}\label{delphilest1}
\inf_{z\in \DD}|\del \phi(z)|\geq K_2
\end{equation}
for some constant $ K_2=K_2(\kappa,\Gamma)$.
But $ \del_z f_{\mu}= ik f_{\mu}\del_z \phi$. Hence, using (\ref{delphilest1}) and the lower bound for $f_{\mu}$ for $z\in \DD$ as mentioned above, we get $ \inf_{z\in\DD}|J_{f_\mu}(z,k)|\geq C_2$  for some constant $C_2=C_2(\kappa,\Gamma,|k|)$.   \hfill$\Box$


\section{Stability of the complex geometric optics solutions }

In this section we consider  two conductivities $\gamma_1,\gamma_2\in C^\varpi(\DD)$ that are 1 near $\partial\DD$ so that  $\mu_j=(1-\gamma_j)/(1+\gamma_j)$ has compact support in $\DD$ and we can
apply the results of the previous two sections; this restriction on $\gamma_j$ will be removed in the next section. For fixed $k\in\CC$ we want to study the stability of the geometric optics solutions
 \eqref{geom-optics-condn1} but, as in \cite{AP} and \cite{BFR}, we will work with the associated solutions $u_1,u_2$ of $\nabla\cdot\gamma_j\,\nabla u=0$ defined by \eqref{eq:f->u}. Let $\om$ and $\vartheta$ be as in Lemma \ref{decomp-g+h} and let $\rho= \rho_{12}=\|\Lambda_{\gamma_1}-\Lambda_{\gamma_2}\|_{\partial\DD}$ 
 where $\|\cdot\|_{\partial\DD}$ denotes the operator norm $H^{1/2}(\partial\DD)\to H^{-1/2}(\partial\DD)$. The stability function that we seek will be of the form
\begin{equation}\label{def:V}
V_k(\rho):=C_1(k)\,[\vartheta(|\log\rho|/C_2)]^{-a}
\end{equation}
for positive constants $C_1(k),C_2,a$. 
Recalling \eqref{est:vartheta}, we see that $V_k(\rho)\to 0$ as $\rho\to 0$  like a negative power of $\log|\log\rho|$; since we are only interested in $\rho\to 0$,  we henceforth assume $0<\rho<1/2$ so that $|\log\rho|>0$. We want to prove the following.

\begin{theorem} \label{th:u-stability} 
Suppose $\gamma_1,\gamma_2\in C^\varpi(\DD)$ such that $\mu_j=(1-\gamma_j)/(1+\gamma_j)$ has compact support in $\DD$ and satisfies $\|\mu_j\|_\infty\leq\kappa<1$ and $\|\mu_j\|_{C^\varpi}\leq\Gamma$. Then, for every $k\in \CC$ there exists $V_k(\rho)$ of the form \eqref{def:V} with constants $C_1(k)$, $C_2$, and $a$ (depending on $\kappa$ and $\Gamma$) such that 
\begin{equation}\label{eq:u-stability}
\|u_1(\cdot,k)-u_2(\cdot,k)\|_{C^0(\DD)}\leq V_k(\rho) \quad\hbox{as}\ \rho\to 0. 
\end{equation}
\end{theorem}
\noindent

We can write
\begin{subequations}
\begin{equation}
u_j(z,k)=e^{ik(z+\e_j(z,k))}
\end{equation}
where (by Theorem \ref{th:estimatefornonlinearBeltrami}) we have $C_*,a>0$ such that
\begin{equation}\label{remainder_k->infty}
|\e_j(z,k)|\leq C_*\,[\vartheta(|k|)]^{-a}
\ \ \hbox{for all $z\in\CC$ and all $|k|> 2$}.
\end{equation}
\end{subequations}
As in \cite{AP}, let us introduce
\begin{subequations}
\begin{equation}\label{def:g(z,w)}
g(z,w,k):= i (z-w)+k\,\e_{z,w}(k),
\end{equation}
where
\begin{equation}\label{def:epsilon_(z,w)}
\e_{z,w}(k):= i(\e_1(z,k)-\e_2(w,k)).
\end{equation}
\end{subequations}
We claim the following is true:
\begin{proposition}\label{pr:g=0=>est}  
For  $C_*$ and $a$ as in \eqref{remainder_k->infty}, there is a constant $C_1>0$ so that $g(z,w,k)=0$ for some $k\not=0$ implies $|z-w|\leq C_1\,[\vartheta(|\log\rho|/4C_*)]^{-a}$.
\end{proposition}
\noindent
This proposition was obtained as Prop.\ 5.3 in \cite{BFR}  for the H\"older case $\om(r)=r^\gamma$, $\vartheta(r)\approx c\,r^{2\gamma}$, $0<\gamma<1$. The proof of Proposition \ref{pr:g=0=>est} follows the same outline; but, 
for completeness, we explain this in the Appendix, including some  details that were
 missing in  \cite{BFR}.
Now let us use Proposition \ref{pr:g=0=>est} to prove our theorem.

\medskip\noindent
{\bf Proof of Theorem \ref{th:u-stability}.}
For $k=0$, $u_j(z,0)=1$, so the left hand side of \eqref{eq:u-stability} is zero. Hence
let us fix $k\not=0$ and pick $z\in\DD$. Using the fact that $\delta_1(\cdot,k)$ is onto $\CC$ (cf.\ Prop.\ 5.2 in \cite{AP}),  there is a $w\in\CC$ such that $\delta_1(w,k)=\delta_2(z,k)$ and hence $g(z,w,k)=0$.
 Let $U$ be a bounded, open set containing both $z$ and $w$.
Then by Theorem \ref{th:f-regularity} in the previous section, we  know that $u_1(z,k)$ is $C^1$ on $\overline U$, so
\[
|u_1(z,k)-u_2(z,k)|=|u_1(z,k)-u_1(w,k)| 
\leq C(k)\,|w-z|.
\]
Proposition \ref{pr:g=0=>est} shows that $|w-z|\leq C_1\,[\vartheta(|\log\rho|/4C_*)]^{-a}$,
so we have  \eqref{eq:u-stability}. \ $\Box$


\section{Proof of the Main Theorem}\label{proofMain}

Now we return to a bounded Lipschitz domain $U$ which we may assume satisfies $\overline{U}\subset\DD$.
For $\gamma_1,\gamma_2\in C^\varpi(\overline{U})$, we want to be able to assume that $\gamma_j\in C^\varpi(\DD)$ with $\gamma_j=1$ near $\partial\DD$ so we can apply our results from Sections 3-5.
This can be achieved using the Whitney extension (cf.\ \cite{S}). As in the Introduction, let $\|\cdot\|_{\partial U}$ denote the norm of an operator $H^{1/2}(\partial U)\to H^{-1/2}(\partial U)$, but now also let  $\|\cdot\|_{\partial\DD}$ denote 
 the norm of an operator $H^{1/2}(\partial \DD)\to H^{-1/2}(\partial \DD)$. The following is the analogue of Theorem 6.2 in \cite{BFR}. 
 \begin{theorem}\label{WhtineyExt} 
 Let $U$ be a Lipschitz domain satisfying $\overline{U}\subset\DD$ and 
 $\gamma_1,\gamma_2\in C^\varpi(\overline{U})$ satisfying $\|\gamma_j\|_{C^\varpi(\overline{U})}\leq \Gamma$ and $\|\gamma_j\|_\infty\leq\kappa <1$. There exists a constant $C=C(\kappa,U)$ and extensions $\widetilde\gamma_1,\widetilde\gamma_2 $ to $\DD$ such that  {\rm supp}$(\widetilde\gamma_j-1)\subset\DD$, $\|\widetilde\gamma_j\|_{C^\varpi(\overline{\DD})}\leq C\,\Gamma$, and
 \begin{equation}\label{est:DtoNcomparison}
 \|\Lambda_{\widetilde\gamma_1}-\Lambda_{\widetilde\gamma_2}\|_{\partial\DD}\leq 
 C\, \|\Lambda_{\gamma_1}-\Lambda_{\gamma_2}\|_{\partial U}.
 \end{equation}
 \end{theorem}
 \noindent
The proof of this theorem follows the same steps as in \cite{BFR} so we will not  discuss it.

\medskip\noindent
{\bf Proof of Main Theorem.} Given $\gamma_1,\gamma_2\in C^\varpi(\overline{U})$, we may use Theorem \ref{WhtineyExt} to consider $\gamma_1,\gamma_2\in C^\varpi(\overline{\DD})$ such that
$\mu_j=(1-\gamma_j)/(1+\gamma_j)$ has compact support in $\DD$.  For $k\in\CC$, let $f_j(z,k)=f_{\mu_j}(z,k)$ be the  complex geometric optics solution  of the associated Beltrami equation $\bar\del f=\mu\,\overline{\del f}$ that was discussed in Section \ref{Estimates-k}.
Let
 \begin{equation}\label{def:U}
 {\mathcal F}(k)= {\mathcal F}(\cdot,k)=f_1(\cdot,k)-f_2(\cdot,k).
 \end{equation}
By  (\ref{eq:regularity-f}) we know that
\begin{subequations}
\begin{equation}\label{est:C1sigma}
\| {\mathcal F}(k)\|_{C^{1,\sigma}(\DD)} \leq C(|k|).
\end{equation}
On the other hand, by Theorem \ref{th:u-stability} there exists
$V_k(\rho)$ of the form \eqref{def:V} so that 
\begin{equation}\label{eq:f-stability}
\| {\mathcal F}(k)\|_{C^0(\DD)}\leq V_k(\rho), 
\end{equation}
\end{subequations}
where $\rho=\rho_{\DD}= \|\Lambda_{\widetilde\gamma_1}-\Lambda_{\widetilde\gamma_2}\|_{\partial\DD}$.
We need to interpolate between \eqref{est:C1sigma} and \eqref{eq:f-stability} to show that  $\|{\mathcal F}(k)\|_{C^1(\DD)}\to 0$ as $\rho\to 0$. In fact, we only need this for one nonzero value of $k$, so 
 let us fix $k=1$ and indicate the $\rho$-dependence by ${\mathcal F}_\rho$. Then \eqref{est:C1sigma} and \eqref{eq:f-stability} become
 \begin{subequations}
 \begin{equation}\label{est:C1sigma-1}
\| {\mathcal F}_\rho\|_{C^{1,\sigma}(\DD)} \leq C(1),
\end{equation}
 \begin{equation}\label{eq:f-stability-1}
\| {\mathcal F}_\rho\|_{C^0(\DD)}\leq V_1(\rho).
\end{equation}
 \end{subequations}
 We want to interpolate to show the spatial derivatives $D{\mathcal F}_\rho$ satisfy
\begin{equation}\label{est:C1}
\|D{\mathcal F}_\rho\|_{C^{0}(\DD)} \leq  V^*(\rho),
\end{equation}
where $V^*(\rho)$ is a nondecreasing positive function  satisfying $V^*(\rho)\to 0$ as $\rho\to 0$. 

Proving \eqref{est:C1} is somewhat technical and uses a  proposition that we have proved in the Appendix. Multiplying ${\mathcal F}_\rho$ by a smooth cutoff function which is 1 on $\DD$, we can assume ${\mathcal F}_\rho\in C^{1,\sigma}_0(\CC)$. Applying Proposition \ref{pr:Sperner} in the Appendix, we obtain
\begin{equation}\label{est:F-interpolated}
\|D{\mathcal F}_\rho\|_{C^0}\leq 2\,\sigma\left(\zeta^{-1}\left(\|{\mathcal F}_\rho\|_{C^0}/[D{\mathcal F}_\rho]_\sigma\right)\right) [D{\mathcal F}_\rho]_\sigma,
\end{equation}
where $\zeta(r):=r\sigma(r)$ is strictly increasing $[0,\infty)\to [0,\infty)$ and surjective so its inverse $\zeta^{-1}$ is well-defined and also strictly increasing; the notation 
$[f]_\sigma $ is defined in \eqref{def:[f]_sigma}.
We need to show the right hand side tends to zero as $\rho\to 0$ in order to conclude \eqref{est:C1}.
We know by \eqref{eq:f-stability-1} that $\|{\mathcal F}_\rho\|_{C^0}\to 0$ as $\rho\to 0$, so $\sigma\left(\zeta^{-1}\left(\|{\mathcal F}_\rho\|_{C^0}/[D{\mathcal F}_\rho]_\sigma\right)\right)\to 0$
as $\rho\to 0$, provided $[D{\mathcal F}_\rho]_\sigma\geq\e>0$ as $\rho\to 0$. However, if $[D{\mathcal F}_\rho]_\sigma\to 0$ as $\rho\to 0$, we still know the right hand side of \eqref{est:F-interpolated} tends to zero because $\sigma$ is bounded on $[0,\infty)$. So \eqref{est:C1} holds.

Now we can use $\mu=\del_{\overline{z}}f_\mu/\overline{\partial_z f_\mu}$ and the lower bound
$\inf_{\DD}|\partial_z f|\geq m$ provided by \eqref{delphilest1} to estimate
\[
\begin{aligned}
\|\gamma_1-\gamma_2\|_{C^0(U)} & \leq \|\gamma_1-\gamma_2\|_{C^0(\DD)} \\
& \leq  \frac{4}{1-\kappa^2}\,\|\mu_1-\mu_2\|_{C^0(\DD)} \\
& \leq  \frac{4}{(1-\kappa^2)m}\,\| Df_1-Df_2\|_{C^0(\DD)}.
\end{aligned}
\]
But finally we  use \eqref{est:C1}, \eqref{est:DtoNcomparison}, and the fact that $V^*$ is nondecreasing to conclude 
\[
\|\gamma_1-\gamma_2\|_{C^0(U)}\leq V^*(\rho_{\DD})
\leq V^*(C\,\rho_{U})
\]
where $\rho_U= \|\Lambda_{\gamma_1}-\Lambda_{\gamma_2}\|_{\partial U}$ and similarly for $\rho_\DD$. This completes the proof. \ $\Box$


\appendix
\section{Additional Lemmas and Proofs}

\begin{lemma}\label{technical-lemma} 
\[
\int_{-\pi}^\pi (1-\cos[2\pi s \cos\theta]) d\theta \geq c \quad \hbox{for all}\ s\geq 1.
\]
\end{lemma}
\noindent
{\bf Proof.}
Since $\cos\theta$ is an even function, it suffices to prove
\begin{equation}\label{eq:F(s)=int(dtheta)}
F(s):=\int_0^\pi \cos[2\pi s\cos\theta]\,d\theta \leq \pi -\e \quad\hbox{for all}\ s\geq 1.
\end{equation}
To do this, let us introduce a change of variables $x=s\cos\theta$, so $x$ ranges from $s$ to $-s$ as $\theta$ ranges from
$0$ to $\pi$, and $d\theta = -\frac{dx}{\sqrt{s^2-x^2}}$.
This means that we can write
\begin{equation}\label{eq:F(s)=int(dx)}
F(s)=\int_{-s}^s \cos[2\pi x]\frac{dx}{\sqrt{s^2-x^2}}=2\int_{0}^s \cos[2\pi x]\frac{dx}{\sqrt{s^2-x^2}}.
\end{equation}
If we estimate $F(s)$ using $|\cos[2\pi x]|\leq 1$, then we obtain $ F(s)=\pi,$
which is not good enough. So we need to make use of values of $x$ for which $\cos[2\pi x]$ is negative.

For fixed $s\geq 1$, we note that $f_s(x)=1/\sqrt{s^2-x^2}$ is increasing in $x$
for $0<x<s$.
Consequently, although $\cos[2\pi x]$ is positive for $0<x<1/4$, we may conclude that
\[
\int_0^{1/4}  \cos[2\pi x]\frac{dx}{\sqrt{s^2-x^2}} + \int_{1/4}^{1/2}  \cos[2\pi x]\frac{dx}{\sqrt{s^2-x^2}}<0.
\]
Trivially, we then conclude that
$\int_0^{3/4}  \cos[2\pi x]\frac{dx}{\sqrt{s^2-x^2}}<0.$
For the same reason, if $s>7/4$, we have
$\int_{3/4}^{7/4}  \cos[2\pi x]\frac{dx}{\sqrt{s^2-x^2}}<0$
and we may add the integrals together to conclude the negativity of the integral over $(0,7/4)$. 
Generalizing this, we conclude that 
\[
\int_{0}^{[s]-1/4}  \cos[2\pi x]\frac{dx}{\sqrt{s^2-x^2}}<0,
\]
where $s\geq 1$ and $[s]$ denotes the greatest integer less than or equal to $s$. Thus
\begin{equation}\label{eq:F(s)-est}
F(s)<2\int_{[s]-\frac{1}{4}}^s \frac{dx}{\sqrt{s^2-x^2}}=\pi -2\sin^{-1}\frac{[s]-\frac{1}{4}}{s}.
\end{equation}
But it is easy to see that
$\frac{[s]-\frac{1}{4}}{s}\geq \eta >0 \quad\hbox{for all}\ s\geq 1, $
and $\sin^{-1}(t)$ is a positive and increasing function for $0<t<1$, so from (\ref{eq:F(s)-est}) we conclude that $F(s)\leq \pi -\e$ as desired. $\Box$

\medskip
As in the proof of Proposition \ref{pr:estimateforlinearBeltrami}, let 
$
K_z(y)=\frac{1}{\pi} \chi(|y|) (z-y)^{-1},
$
where $\chi(r)$ is a smooth cut-off function satisfying $\chi(r)=1$ for $0\leq r\leq 1$ and $\chi(r)=0$ for $r\geq 3/2$.
Since $K_z\in L^1(\RR^2)$ we know that $\widehat{K_z}(\xi)$ is bounded for $\xi\in\RR^2$ and even 
$\widehat{K_z}(\xi)\to 0$ as $|\xi|\to \infty$ by Riemann-Lebesgue. But we need 
 more precise decay as $|\xi|\to \infty$.
 
\begin{lemma}  
The Fourier transform $\widehat{K_z}(\xi)$ satisfies the estimate \eqref{est:hat-K}.
\end{lemma}
\noindent
{\bf Proof.} For $|z|>3/2$, $K_z(y)$ is a smooth function of $y\in\RR^2$, so $\widehat{K_z}(\xi)$ decays rapidly as $|\xi|\to\infty$, uniformly for $|z|\geq 2$. Se we restrict our attention to $|z|\leq 2$.
It suffices to consider $\xi=(s,0)$ for $s>2$ and use the equivalent definition of the Fourier transform,
$\widehat f(\xi)=\int e^{iy\cdot\xi}f(y)\,dy$.
In the following, we use polar coordinates $re^{i\theta}$ for $u:=y-z$ and observe that $|y|<2$, $|z|\geq 2$ imply $|u|<3$:
\[
-\pi \widehat{K_z}(\xi)=\int_{|y|<2} \frac{e^{iy\cdot\xi}\,\chi(|y|)}{y-z}\,dy
= e^{iz\cdot\xi}\int_0^{2\pi} e^{-i\theta} \int_0^4 e^{isr\cos\theta}\chi(re^{i\theta}+z)\,dr \,d\theta.
\]
Now let us integrate by parts:
\[
 \int_0^4 e^{isr\cos\theta}\chi(re^{i\theta}+z)\,dr=
 -\frac{1}{is\cos\theta}\int_0^4\left(e^{isr\cos\theta}-1\right)\,\frac{d}{dr}\chi(r\,e^{i\theta}+z)\,dr,
\]
where we observe that $\chi(re^{i\theta}+z)$ is constant for $0<r<\e$.
For $\theta\in (0,\pi/2)$ let us make the substitution $t=\cos\theta$, $dt=-\sin\theta\,d\theta$, to find
\[
\begin{aligned}
 \left|  \int_0^{\pi/2} \int_0^4 e^{isr\cos\theta}\,\chi(re^{i\theta}+z)\,dr\,d\theta \right|
& \leq \frac{1}{s}\int_0^1\int_0^4 \frac{|e^{isrt}-1|}{t\sqrt{1-t^2}}\,|f(r)|\,dr\,dt \\
& \leq \frac{C}{s} \int_0^4\int_0^1 \frac{|\cos(srt)-1|+|\sin(srt)|}{t\sqrt{1-t^2}} \,dt\,dr,
\end{aligned}
\]
where $f(r)=f_{\theta,z}(r)= \frac{d}{dr}\chi(r\,e^{i\theta}+z)$ and the constant $C$ depends on the maximum of $f$. Let us focus on the integral involving $\sin(srt)$. If $sr<1$ then
\[
\int_0^1  \frac{|\sin(srt)|}{t\sqrt{1-t^2}}\,dt\leq \int_0^1 \frac{dt}{\sqrt{1-t}}=1.
\]
If $sr>1$ then 
\[
\int_{0}^1 \frac{|\sin(srt)|}{t\sqrt{1-t^2}}\,dt
\leq \int_0^{1/{sr}} \frac{sr \,dt}{\sqrt{1-t^2}} + 
\int_{1/{sr}}^1 \frac{dt}{t\,\sqrt{1-t}}.
\]
We can evaluate the first of these integrals and then estimate as $sr\to\infty$:
\[
\int_0^{1/{sr}} \frac{sr \,dt}{\sqrt{1-t^2}} = sr\sin^{-1}\frac{1}{sr}
\approx 1.
\]
For the second integral, we can use an integral table and $\sqrt{1-a}\approx 1-\frac{a}{2}$ as $a\to 0$:
\[
\begin{aligned}
 \int_{1/{sr}}^1 \frac{dt}{t\,\sqrt{1-t}}
&=  \log \left| \frac{\sqrt{1-1/{sr}}+1}{\sqrt{1-1/{sr}}-1}\right| \\
&= \log sr +O(1)\leq \log s +\log|r|+ C \quad\hbox{for}\ s>2,\ 0<r<4.
\end{aligned}
\]
We conclude that
\[
\int_{0}^1 \frac{|\sin(srt)|}{t\sqrt{1-t^2}}\,dt
\leq C\,(\log s + \log r) \quad \hbox{for} \ s>2,\ 0<r<4.
\]
We can similarly show
\[
 \int_{0}^1 \frac{|\cos(st)-1|}{|t|\sqrt{1-t^2}}\,dt \leq C\,(\log s + \log r) \quad \hbox{for} \ s>2,\ 0<r<4,
\] 
so we conclude 
\[
\left|  \int_0^{\pi/2} \int_0^4 e^{isr\cos\theta}\,\chi(re^{i\theta}+z)\,dr\,d\theta \right|
\leq C\, \frac{\log s }{s} \quad\hbox{for} \ s>2.
\]
The substitution $t=\cos s$ can be used again for $\theta\in (\pi/2,\pi)$, $(\pi,3\pi/2)$ and $(3\pi/2,2\pi)$, so we 
can put these all together to obtain the desired estimate:
\[
\pi |\widehat{K_z}(\xi)| \leq C\,\frac{\log|\xi|}{|\xi|}\quad\hbox{for}\ |\xi|>2. \quad \Box
\]

\medskip  
Now we begin the preparations to prove Proposition \ref{pr:g=0=>est} on the
complex gemoetric optics solutions $u_j(z,k)=\exp[ik(z+\e_j(z,k)]$ which satisfy
$|\e_j(z,k)|\leq C_*[\vartheta(|k|)]^{-a}$ for all $z\in\CC$ and $|k|\geq 2$.
To study $g(z,w,k)=i(z-w)+k\e_{z,w}(k)$ where $\e_{z,w}(k)=i(\e_1(z,k)-\e_2(z,k))$, we need to treat its behavior for large $|k|$ differently from small $|k|$; but what is ``large'' and what is  ``small" depends on 
\[
\lambda:=z-w.
\]
As in \cite{BFR}, we want to define a function $R:\CC\to \RR$ so that for $|k|\geq R(\lambda)$ we have 
 $|\e_{z,w}(k)|\leq |\lambda|/2$ and hence $g(z,w,k)\not=0$. In fact, since $\vartheta:[1,\infty)\to[0,\infty)$ is strictly increasing,
 let us denote its inverse by $\vartheta^{-1}:[0,\infty)\to [1,\infty)$. If we use the above constants $C_*$ and $a$ to define
\begin{equation}
R(\lambda):= \vartheta^{-1}\left(\left| \frac{\lambda}{4C_*}\right|^{-1/a}\right),
\end{equation}
then $|k|\geq R(\lambda)$ indeed implies
\[
|\e_{z,w}(k)|\leq 2\,C_*\,\left[\vartheta\left(|k|\right)\right]^{-a}\leq 2\,C_*\,\left[\vartheta\left(R(\lambda)\right)\right]^{-a}
=\frac{|\lambda|}{2}.
\]
The proof of Proposition \ref{pr:g=0=>est} is then reduced to finding a constant $C_1$ so that if
\begin{equation}\label{lambda-lowerbound}
|\lambda |> C_1\,[\vartheta(|\log\rho|/4C_*)]^{-a},
\end{equation}
then the only zero of $g(z,w,k)$ in the set $|k|\leq R(\lambda)$ is at $k=0$. 
The following is a simple relationship between $\lambda$ and $\rho$ that is useful in
subsequent proofs. 

\begin{lemma} 
There is a constant $C_1$ such that if $\lambda$ satisfies \eqref{lambda-lowerbound}, then
\begin{equation}\label{rho-upperbound}
\rho <|\lambda|\,e^{-C_*R(\lambda)}.
\end{equation}
\end{lemma}
\noindent
{\bf Proof.}
Note that \eqref{lambda-lowerbound} implies
$
 |\lambda/C_1|^{-1/a} < \vartheta(|\log\rho|/4C_*),
$
and the strict monotonicity of $\vartheta^{-1}$ implies
$
\vartheta^{-1}(|\lambda/C_1|^{-1/a})< |\log\rho|/4C_*.
$
Now, provided $C_1\geq 4C_*$, we have $|\lambda/4C_*|^{-1/a}\leq |\lambda/C_1|^{-1/a}$, so
by the monotonicity of $\vartheta^{-1}$ we have
\[
R(\lambda)=\vartheta^{-1}\left(\left|\frac{\lambda}{4C_*}\right|^{-1/a}\right)
\leq \vartheta^{-1}\left(\left|\frac{\lambda}{C_1}\right|^{-1/a}\right)<\frac{|\log\rho|}{4C_*}.
\]
Consequently, we have $e^{-C_*R(\lambda)}>\rho^{1/4}.$
Thus to obtain \eqref{rho-upperbound} it suffices to show $\rho^{3/4}<|\lambda|$. Using \eqref{lambda-lowerbound} again, we see that it suffices that $C_1$ is an upper bound for
\begin{equation}
f(\rho)=\rho^{3/4}  \left[\vartheta(|\log\rho|/4C_*)\right]^{a}
  \quad\hbox{for}\ 0<\rho<e^{-4C_*}.
\end{equation}
But we know from \eqref{est:vartheta} that $\vartheta(|\log\rho|/4C_*)$ grows only like a power of 
$\log |\log\rho |$ as $\rho\to 0$, so $f(\rho)\to 0$ as $\rho\to 0$. Thus such a $C_1$ may be found.
$\Box$

\medskip
For fixed $z,w$, the function $g$ satisfies a  $\bar\partial$-equation in the variable $k$. Since this does not involve the regularity of $\mu$, we may import results from the H\"older case. The following appears as  Lemma 5.4 in  \cite{BFR}:
\begin{lemma}\label{le:d-bar(g)} 
For fixed $z,w$, the function $g$ as in \eqref{def:g(z,w)} satisfies
\begin{equation}\label{d-bar(g)}
\partial_{\overline k} g=\sigma g + E,
\end{equation}
where $\sigma=\sigma_{z,w}$ and $E=E_{z,w}$ satisfy
\begin{equation}\label{sigma,E-bounds}
|\sigma(k)|\leq 2, \quad
|E(k)|\leq \rho\,e^{c_1(1+|k|)}, \quad
|D E(k)|\leq e^{c_1(1+|k|)},
\end{equation}
for some constant $c_1=c_1(\kappa)>0$.
\end{lemma}

We want to obtain conclusions about the behavior of $g$ from the fact that it satisfies \eqref{d-bar(g)}. This requires inverting the operator $\partial_{\overline{k}}$, but we do not have sufficient decay at infinity to directly apply the Cauchy transform, so we need to multiply by a cut-off function. For the moment, let us ignore $\lambda$. We fix $R\geq 2$ and consider a cut-off function $\chi_R\in C_0^\infty(\DD_{2R})$ with $\chi_R(k)=1$ for $|k|\leq R$. 
Then, for the functions $\sigma$ and $E$ in Lemma \ref{le:d-bar(g)}, let us introduce
\begin{equation}\label{def:eta,S}
\eta_R(k):= P(\sigma\chi_R) \quad\hbox{and}\quad S_R(k):=P(e^{-\eta_R}E\chi_R).
\end{equation}
Here, $P$ denotes the Cauchy transform (in the variable $k$), so
\[
\partial_{\overline k} \,\eta_R = \sigma \chi_R \quad\hbox{and}\quad \partial_{\overline k}\,S_R= e^{-\eta_R}E\chi_R.
\]
The functions $\eta_R$ and $S_R$ have the following global estimates: cf.\ Lemma 5.5 in \cite{BFR}.

\begin{lemma}\label{etaS} 
For fixed $z,w$, there is a constant $c_2=c_2(\kappa)$  such that
\begin{equation}\label{eta,S-estimates}
\|\eta_R\|_{L^\infty(\CC)}  \leq c_2 R \quad\hbox{and}\quad \| S_R\|_{L^\infty(\CC)} \leq \rho\,  e^{c_2R}.
\end{equation}
In fact, for any $0<\theta<1$, $c_2$ may be chosen so that
\begin{equation}\label{gradS-estimate}
\|\nabla S_R\|_{L^\infty(\CC)} \leq \rho^\theta\,e^{c_2\,R}.
\end{equation}
\end{lemma}

Now let
$
\widetilde S_R(k)=S_R(k)-S_R(0)
$
so that $\widetilde S_R(0)=0$. Then $\widetilde S_R$  satisfies \eqref{gradS-estimate} and $\|\widetilde S\|_{L^\infty}\leq 2\rho e^{c_2R}$. Let us define
\begin{equation}\label{def:F}
F(k)\equiv F(z,w,k):=e^{-\eta_R(k)}\,g(z,w,k)-\widetilde S_R(k).
\end{equation}
A straightforward calculation shows $\partial_{\overline k}F(z,w,k)=0$ for $|k|\leq R$, so 
for fixed $z,w$ the function $F$ is analytic for $k\in\DD_R$. By construction, $F(0)=0$ and the following result shows that this is the only zero in $\DD_R$ when $\lambda$ is sufficiently large.
\begin{lemma} \label{F-zeros} 
There is a constant $C_1$ such that for $\lambda$ satisfying \eqref{lambda-lowerbound} the function $F(z,w,k)$ has a unique zero at $k=0$ in the set $|k|\leq R(\lambda)$.
\end{lemma}
\noindent
This appears as  Proposition 5.6 in \cite{BFR}, which is proved by showing that $F(k)$ is homotopic to $k$ through nonvanishing functions on $|k|=R(\lambda)$ and uses the estimate \eqref{rho-upperbound}; thus it may be repeated in our case. 

\medskip
We shall need some additional properties of $F$.
\begin{lemma} \label{F=approx} 
For $\lambda$ satisfying \eqref{lambda-lowerbound} with $C_1$ given in Lemma \ref{F-zeros}, we can write
\[
F(z,w,k)=  \lambda  k \, e^{\nu(k)} \quad\hbox{for}\ |k|\leq R(\lambda),
\]
where $\nu(k)$ is analytic and satisfies
$
|\nu(k)|\leq {c_2 R(\lambda)}
$
with $c_2$ from Lemma \ref{etaS}.
\end{lemma}
\noindent
This appears as Lemma 5.7 in \cite{BFR}, which is proved using the analyticity of $F(k)/k$ in
 $|k|<R(\lambda)$; the estimate \eqref{rho-upperbound}  is used with the maximum principle  and may be repeated in our case.
The next two results follow from Lemma \ref{F=approx} exactly as in  \cite{BFR}.

\begin{corollary}\label{co:F-inverse} 
For $\lambda$ satisfying \eqref{lambda-lowerbound} with $C_1$ given in Lemma \ref{F-zeros} and any $\delta>0$, we have
\[
F^{-1}(\DD_\delta)\subset \DD_{\delta_1} \quad\hbox{where}\ \delta_1=\delta\, e^{c_2 R(\lambda)}/|\lambda|,
\]
with $c_2$  from Lemma \ref{etaS}.
\end{corollary}

\begin{corollary}\label{co:F'-lowerbound} 
For $\lambda$ satisfying \eqref{lambda-lowerbound} with $C_1$ given in Lemma \ref{F-zeros}, there exists $d>0$ so that
\[
\inf_{|k|<d} |F'(k)| >\frac{1}{2} |\lambda| e^{-c_2 R(\lambda)}.
\]
\end{corollary}

\medskip
Now in order to reach our desired conclusion that, for $\lambda=z-w$ satisfying \eqref{lambda-lowerbound}, $g(z,w,k)$ only vanishes at $k=0$, it suffices to show that the function
\begin{equation}
H_{z,w}(k)=e^{-\eta(k)}g(z,w,k)=F_{z,w}(k)+\widetilde S_{z,w}(k)
\end{equation} 
has a unique zero at $k=0$ in the set $|k|\leq R(\lambda)$. Note that $H$ is not analytic in $|k|\leq R(\lambda)$,
so we cannot use the principle of the argument as we did in the proof of Lemma \ref{F-zeros};
instead we shall apply degree theory to $H$. For this we need to know more about the zeros of $H$: for given values of $z,w$ let
\[
Z(H_{z,w})=\{k\in\CC:H(z,w,k)=0\}.
\]
The following two Lemmas and Proof of Proposition \ref{pr:g=0=>est} follow the ideas in \cite{BFR}.
\begin{lemma}\label{le:Z(H)inDisk} 
There exists $C_1$ such that for $\lambda$ satisfying \eqref{lambda-lowerbound}  we have 
\[
Z(H)\subset \DD_d,
\]
where $d$ is given in Corollary \ref{co:F'-lowerbound}.
\end{lemma}
\noindent
{\bf Proof.} For $k\in Z(H)$, $F(k)=-S(k)$. But by Lemma \ref{etaS}, we have $\|S\|_\infty \leq \rho\,e^{c_2R(\lambda)}$. So by Corollary \ref{co:F-inverse}, $|k|<\rho\, e^{2c_2 R(\lambda)}/|\lambda|$. Thus if we have chosen $\lambda$ so that 
$\rho\,e^{2c_2 R(\lambda)}/|\lambda|<d$, then we will have $Z(H)\subset \DD_d$. But, recalling that $R(\lambda)\to 1$ as $|\lambda|\to\infty$, this can be arranged by requiring 
$|\lambda |> C_1\,\vartheta(|\log\rho|/4C_*)^{-a}$ for $C_1$ sufficiently large.  \hfill$\Box$

\medskip
We also need to know about the Jacobian determinant of $H$, which can be expressed (cf.\ \cite{AIM}) as
\begin{equation}\label{det-DH}
\hbox{det} \, DH=|H_k|^2-|H_{\bar k}|^2.
\end{equation}
\begin{lemma} \label{le:sign-det-DH} 
There is a constant $C_1$ such that for $\lambda$ satisfying \eqref{lambda-lowerbound}  we have \hbox{\rm det}\,$DH(k)>0$ for all 
$|k|<d$, where $d$ is as in Corollary \ref{co:F'-lowerbound}.
\end{lemma} 
\noindent
{\bf Proof.} Since $H=F+S$ with $F$ analytic for $|k|\leq R(\lambda)$, we have $H_k=F_k+S_k=F'+S_k$ and $H_{\bar k}=S_{\bar k}$.
Using $-2\,\hbox{Re}(F'S_k)\leq 2|F'|\,|S_k|\leq \frac{1}{2}|F'|^2+2|S_k|^2$, we can easily show
\[
\hbox{det} \,DH \geq \frac{1}{2}|F'|^2 - |DS|^2.
\]
Now assuming $|k|<d$ so that we can use Corollary \ref{co:F'-lowerbound} and using \ref{gradS-estimate} we have 
\[
 \frac{1}{2}|F'|^2 - |DS|^2 > \frac{1}{2} |\lambda| e^{-c_2 R(\lambda)} - \rho^{2\theta}\,e^{2\,c_2\,R(\lambda)}.
\]
Thus we can prevent det\,$DH$ from vanishing by choosing $|\lambda|$ large enough that 
\[
 |\lambda | e^{-3c_2\,R(\lambda)}>2\,\rho^{2\theta}.
\]
Recalling that $R(\lambda)\to 1$ as $|\lambda|\to\infty$, we see that this can be achieved by taking $|\lambda |> C_1\,[\vartheta(|\log\rho|/4C_*)]^{-a}$ with $C_1$ sufficiently large. \hfill$\Box$

\medskip
Now we are finally ready to give the proof of Proposition \ref{pr:g=0=>est}.

\medskip\noindent
{\bf Proof of Proposition \ref{pr:g=0=>est}:}
To begin with, let $\C_R=\partial\DD_R$ denote the circle of radius $R$. For $\lambda$ satisfying \eqref{lambda-lowerbound} we know that $g=i\lambda k+\e(k)k$ is homotopic to $i\lambda k$, and hence to $k$, through nonvanishing  functions on $\C_R$. Since $H=e^{-\eta}g$ is homotopic to $g$ through nonvanishing functions on $\C_R$, we have
\[
\hbox{deg}(H,\DD_{R(\lambda)},0)=\hbox{deg}(g,\DD_{R(\lambda)},0)=1.
\]
On the other hand, by the degree formula we have
\[
\hbox{deg}(H,\DD_{R(\lambda)},0)=\sum_{k_i\in Z(H)} \hbox{sign det} DH(k_i).
\]
However, we know by Lemmas \ref{le:Z(H)inDisk} and \ref{le:sign-det-DH} that sign\,det\,$DH(k_i)$=+1 for all $k_i\in Z(H)$. So, in order to have 
$\hbox{deg}(H,\DD_{R(\lambda)},0)=1$, we must have only one zero, namely at $k=0$.  \hfill$\Box$

\medskip
The following lemma and proposition are used in the proof of the main theorem.
Let $\sigma(r)$ be a modulus of continuity; for example, using \eqref{defsigma} as in Section \ref{sec:regularity}. (Note that we {\em not} require $\sigma$ to satisfy the Dini condition at $r=0$.)
For $f\in C_0^{\sigma}(\RR^2)$,  let us introduce
\begin{equation}\label{def:[f]_sigma}
[f]_\sigma = \sup_{x\not= y}\frac{|f(x)-f(y)|}{\sigma(|x-y|)}.
\end{equation}
The following result and its proof are taken from \cite{Sp} (cf.\  (10.3) in \cite{Sp}).

\begin{lemma} 
Suppose $f\in C_0^{1,\sigma}(\RR^2)$. Then for any $r>0$ and any $i=1,\dots,n$,
\begin{equation}\label{est:Lemma1}
\|f_{x_i}\|_{C^0}\leq \sigma(r)\left[ f_{x_i} \right]_{\sigma}+\frac{1}{r}\|f\|_{C^0}.
\end{equation}
\end{lemma}

\noindent
{\bf Proof.} For any $y\in \RR^2$ and $r>0$ we can find $y_1,y_2\in \del B(y,r)$ and $\bar y\in B(y,r)$ such that
\[
|f_{x_i}(\bar y)|=\frac{1}{2r}|f(y_1)-f(y_2)|\leq \frac{1}{r}\|f\|_{C^0}.
\]
Thus
\[
\begin{aligned}
|f_{x_i}(y)|&\leq |f_{x_i}(y)-f_{x_i}(\bar y)|  + |f_{x_i}(\bar y)| \\
&\leq \sigma(r) [f]_\sigma + \frac{1}{r}\|f\|_{C^0}.
\end{aligned}
\]
Taking supremum over $y\in\RR^2$ yields \eqref{est:Lemma1}. $\Box$

\medskip
Note that $\zeta(r)=r\sigma(r)$ is strictly increasing, so its inverse function $\zeta^{-1}(r)$ is defined and also strictly increasing. The following is a sort of interpolation inequality.

\begin{proposition}\label{pr:Sperner} 
If  $f\in C_0^{1,\sigma}(\RR^2)$, then 
\begin{equation}\label{est:interpolation}
\|f_{x_i}\|_{C^0} \leq 2\,\sigma\left(\zeta^{-1}\left(\|f\|_{C^0}/[f_{x_i}]_\sigma\right)\right) [f_{x_i}]_\sigma.
\end{equation}
\end{proposition}

\noindent
{\bf Proof.} Since $f$ cannot be identically constant unless it is identically zero, we can assume $\left[ f_{x_i} \right]_{\sigma}\not= 0$. Starting from $r=0$, increase $r$ until $\sigma(r)\left[ f_{x_i} \right]_{\sigma}=\frac{1}{r}\|f\|_{C^0}.$
For this value of $r$, \eqref{est:Lemma1} becomes
\begin{equation}\label{est:r-dependent}
\|f_{x_i}\|_{C^0}\leq 2\,\sigma(r)\left[ f_{x_i} \right]_{\sigma}.
\end{equation}
We need to eliminate $r$ from this estimate. But we know
\[
\zeta(r)=r\sigma(r)=\frac{\|f\|_{C^0}}{\left[ f_{x_i} \right]_{\sigma}}.
\]
So
\[
r=\zeta^{-1}\left(\frac{\|f\|_{C^0}}{\left[ f_{x_i} \right]_{\sigma}}\right).
\]
Plugging this into \eqref{est:r-dependent} yields \eqref{est:interpolation}. $\Box$

\bigskip\noindent
For example, if $\sigma=r^\alpha$ for $\alpha\in (0,1)$, then $\zeta(r)=r\sigma(r)=r^{1+\alpha}$ and $\zeta^{-1}(r)=r^{1/(1+\alpha)}$. Consequently, \eqref{est:interpolation} implies the more familiar interpolation inequality
\[
\|f_{x_i}\|_0\leq 2\,\|f\|_0^\theta \,\|f\|_{1,\alpha}^{1-\theta} \quad\hbox{where $\theta=\alpha/(1+\alpha)$.}
\]


\end{document}